\documentclass[11pt]{article}
\usepackage{latexsym,amsmath,amsfonts}
\DeclareSymbolFont{lettersA}{U}{pxmia}{m}{it}
\DeclareMathSymbol{\piup}{\mathord}{lettersA}{"19}
\normalsize
\begin{document}
\begin{center}
\bigskip
{\LARGE\textbf{On the Degree Sequence\\of Random
Geometric Digraphs}}\\
\bigskip
Yilun Shang\footnote{Department of Mathematics, Shanghai Jiao Tong University, Shanghai 200240, CHINA. email: \texttt{shyl@sjtu.edu.cn}}\\
\end{center}
\bigskip
\begin{abstract}
A random geometric digraph $G_n$ is constructed by taking
$\{X_1,X_2,\cdots X_n\}$ in $\mathbb{R}^2$ independently at random
with a common bounded density function. Each vertex $X_i$ is
assigned at random a sector $S_i$ of central angle $\alpha$ with
inclination $Y_i$, in a circle of radius $r$ (with vertex $X_i$ as
the origin). An arc is present from vertex $X_i$ to $X_j$, if $X_j$
falls in $S_i$. Suppose $k$ is fixed and $\{k_n\}$ is a sequence
with $1\ll k_n\ll n^{1/2}$, as $n\rightarrow\infty$. We prove
central limit theorems for $k-$ and $k_n-$nearest neighbor distance
of out- and in-degrees in $G_n$. We also show that the degree
distribution of this model, which varies with the probability
distribution of the underlying point processes, can be either
homogeneous or inhomogeneous. Our work should provide valuable
insights for alternative mechanisms wrapped in real-world complex
networks.
\end{abstract}

\normalsize
\bigskip

\noindent{\Large\textbf{1. Introduction}}

\smallskip
In random graph theory, degree sequences are among the most
elementary and essential issues. The random geometric graphs
$G(\mathcal{X},r)$ have been well studied in the last decade, see
the monograph \cite{1}, a short overview \cite{21} and references
therein. In order to investigate the typical vertex degree of
$G(\mathcal{X}_n,r_n)$, Penrose(\cite{2}) defined an empirical
process of $k_n$-nearest neighbor distances in $\mathcal{X}_n$, and
showed the weak convergence of the finite-dimensional distributions
of that process, scaled and centered, to a Gaussian limit process.
He further considered the case $k_n=k$ fixed in \cite{1} later.
Given a finite point set $\mathcal{X}\in\mathbb{R}^d$ and given
$x\in\mathcal{X}$, the $k$-nearest neighbor distance means the
distance from $x$ to its $k$-nearest neighbor in $\mathcal{X}$. In
the geometric setting, the $k$-nearest neighbor distance is often a
suitable vehicle to deal with degree-related properties of spatial
point configurations\cite{7}. It is also closely concerned with
$k$-spacing in statistical testing, which has a number of
applications, see the book \cite{22}, and is of interest in its own
right.

In this paper we extend the method of Penrose and establish results
analogous to the ones mentioned above for in-degree and out-degree
of random geometric digraphs. Our result (Theorem 3) shows that the
degree distribution of random geometric digraphs in the
thermodynamic regime can be either homogeneous or inhomogeneous
according to different underlying distributions of point processes.
In particular, the degree distribution is Poisson-like when points
are uniformly scattered, reminiscent of that of Erd\"os-R\'enyi
random graphs, see the classic book \cite{23}(Chap.3); otherwise the
degree distribution is highly skew (or inhomogeneous), similar with
that of many large real-world graphs \cite{25}. We also mention that
the author was recently able to prove the maximum out/in-degrees are
almost determined \cite{30}, and this phenomenon has been discovered
in Erd\"os-R\'enyi random graphs \cite{23}. For more discussions,
see Section 2.1.

Let $\mathcal{X}_n=\{X_1,X_2,\cdots,X_n\}$, $\{X_i\}$ are i.\ i.\
d.\ random variables in $\mathbb{R}^d$ with distribution $F$ having
a specified bounded density function $f$. Let
$\mathcal{P}_n=\{X_1,X_2,\cdots,X_{N_n}\}$, $N_n\sim Poi(n)$. So
$\mathcal{P}_n$ is a Poisson point process with intensity $nf$,
coupled with $\mathcal{X}_n$. Let $\mathcal{H}_\lambda$ be a
homogeneous Poisson process with intensity $\lambda$ on
$\mathbb{R}^d$ and $\|\cdot\|$ be $l^2$ norm on $\mathbb{R}^d$.
Standard random geometric graphs $G(\mathcal{X}_n,r_n)$,
$G(\mathcal{P}_n,r_n)$ are defined as in \cite{1}, that is,
$G(\mathcal{X}_n,r_n)$ (or $G(\mathcal{P}_n,r_n)$) has vertex-set
$\mathcal{X}_n$ (or $\mathcal{P}_n$) and an edge $X_iX_j$
$(i\not=j)$ if $||X_i-X_j||< r_n$. We always assume that
$r_n\rightarrow0$ as $n\rightarrow\infty$. We now define
random geometric digraph models to use in this paper as follows: \\
\textbf{Definition 1.} \itshape $(d=2)$ Let $\alpha\in (0,2\piup]$
be fixed. Let $\mathcal{Y}_n=\{Y_1,Y_2,\cdots,Y_n\}$ be i.i.d.
random variables, taking values in $[0,2\piup)$, with density
function $g$. Associate every point $X_i\in \mathcal{X}_n$ a sector,
which is centered at $X_i$, with radius $r_n$, amplitude $\alpha$
and elevation $Y_i$ with respect to the $x$-axis horizontal
direction anticlockwise. This sector is denoted as $S(X_i,Y_i,r_n)$.
We denote by $G_\alpha(\mathcal{X}_n,\mathcal{Y}_n,r_n)$
(abbreviated as $G_n$) the digraph with vertex set $\mathcal{X}_n$,
and with arc $(X_i,X_j)$, $i\not=j$, present if and only if $X_j\in
S(X_i,Y_i,r_n)$. We can define a Poisson version
$G_{\alpha}(\mathcal{P}_n,\mathcal{Y}_{N_n},r_n)$ ($G'_n$ for short)
similarly.\normalfont

In what follows, we will primarily take
$g=\frac{1}{2\piup}1_{[0,2\piup)}$, that is, $Y_i\sim U[0,2\piup)$.
We will defer the discussion of the case of $d\ge 3$, general
probability density function $g$ and even other norms to Section 6.
Actually, the above model has been first introduced in \cite{3}
under the name ``random scaled sector graph'', with $d=2$, Euclidean
norm and $n$ points uniformly distributed in $[0,1]^2$. This is an
important variant of random geometric graph which has been
revitalized recently in the context of wireless ad hoc networks, and
it is used to analyze the performance of wireless sensor networks
communicating through optical devices or directional antennae, which
are significant in mobile communication\cite{20}. Some basic
properties and graph-theoretic parameters of this model have also
been addressed\cite{3,4,18}, using basically combinatorial
techniques and discretization.

The rest of this paper is organized as follows. Section 2 contains
the statement of main results for $d=2$, $Y_i$ uniformly
distributed. Section 3 discusses the asymptotic results for means
and degree distribution. In Section 4, we give some moments
preparatives for de-Poisson. Section 5 includes the proof of main
theorems. Section 6 is devoted to higher dimension and general
probability density function g.

\bigskip
\bigskip

\noindent{\Large\textbf{2. Statement of main results}}
\smallskip

We will consider two asymptotic regimes. First, take $k_n\equiv k\in
\mathbb{N}$. Second, let $k_n\rightarrow \infty$, and
\begin{equation}
\lim_{n\rightarrow \infty}\frac{k_n}{\sqrt{n}}=0.\label{1}
\end{equation}
Notice that if we want the sequence $\{k_n\}_{n\ge1}$ to converge as
$n$ tends to infinity, then the above two cases are only choices (
and (\ref{1}) is technically needed in the proofs). In the first
regime, define $r_n=r_n(t)$ by $nr_n(t)^2=t$, for $t>0$, and in the
second, define $r_n=r_n(t)$ by $nr_n(t)^2=s(k_n+t\sqrt{k_n})$, for
$s>0$, $t\in \mathbb{R}$. Here we introduce a tunable parameter $t$
to adjust the areas of sectors and $t$ has nothing to do with
``time'', though we will study several random processes with $t$
that evolves. Regulating $t$ allows us to tackle the degree
sequences in fine details. The reason why we choose such $r_n$ is to
ensure a non-degenerate limit, since $nr_n^2$ is a good measure of
average degree, see the appendix A of \cite{5}. We emphasize that
$k_n$ is a crucial parameter which appears in two respects, the
scale on which the degree distribution tails are studied as well as
the scaling for the radius $r_n$.

Before proceeding, we give some notations to ease statement. For
$\lambda >0$, let $\rho_{\lambda}(k):=P(Poi(\lambda)=k)$ and for
$A\subseteq \mathbb{Z}^+$, let $\rho_{\lambda}(A):=P(Poi(\lambda)\in
A)$. For $x\in \mathbb{R}^2$, let $\phi$, $\Phi$ be the density and
distribution function of standard normal variables. Given $x\in
\mathbb{R}^2$, define $B(x,r)$ the disk with center $x$ and radius
$r$, and let $B_n(x,t):=B(x,r_n(t))$, $S_n(x,y,t):=S(x,y,r_n(t))$ in
both limit regimes. Following Penrose \cite{1} we set
$\mathcal{X}^x:=\mathcal{X}\cup\{x\}$, if $\mathcal{X}$ is a finite
set in $\mathbb{R}^2$ and $x\in\mathbb{R}^2$. Denote by
$\#\mathcal{X}$ the number of elements in $\mathcal{X}$ and
$\mathcal{X}(A):=\#(\mathcal{X}\cap A)$ for
$A\subseteq\mathbb{R}^2$.

We will need some further definitions before we can state our main
results. In the rest of the paper $f_{\max}$ will denote the
essential supremum of the probability density function $f$, i.e.
$f_{\max}:=\sup\{u:|\{x:f(x)>u\}|>0\}$. Here and in the rest of the
paper $|\cdot|$ denotes Lebesgue measure. We assume
$f_{\max}<\infty$ throughout the paper. Next, define the level set
when $k_n\rightarrow \infty$ as $L_s:=\{x\in \mathbb{R}^2|s
f(x)=\frac{2}{\alpha}\}$ and let $L_s^+:=\{x\in \mathbb{R}^2|s
f(x)>\frac{2}{\alpha}\}$. We also put a mild restriction on density
function $f$: let $R:=\big\{x\in \mathbb{R}^2|f(x)>0,\
\limsup\limits_{y\longrightarrow
x}\frac{|f(y)-f(x)|}{\|y-x\|}<K\big\}$ with some $K<\infty$, and we
always assume $F(R)=1$. Let $c$, $c'$ be various positive constants,
and the values may change from line to line.

For Borel set $A\subseteq\mathbb{R}^2$, define $\xi_n^{out}(t,A)$,
$\xi_n^{'out}(t,A)$ be the number of vertices in $A$ of out-degrees
at least $k_n$ of $G_n$ and $G'_n$ respectively. More specifically,
\begin{equation}
\xi_n^{out}(t,A)=\sum_{i=1}^n1_{[\mathcal{X}_n(S_n(X_i,Y_i,t))\ge
k_n+1]\bigcap[X_i\in A]} \nonumber
\end{equation}
\begin{equation}
\xi_n^{'out}(t,A)=\sum_{i=1}^{N_n}1_{[\mathcal{P}_n(S_n(X_i,Y_i,t))\ge
k_n+1]\bigcap[X_i\in A]} \nonumber
\end{equation}
Similarly, for in-degree we have,
\begin{equation}
\xi_n^{in}(t,A)=\sum_{i=1}^n1_{[\#\{X_j\in\mathcal{X}_n|X_i\in
S_n(X_j,Y_j,t)\}\ge k_n+1]\bigcap[X_i\in A]} \nonumber
\end{equation}
\begin{equation}
\xi_n^{'in}(t,A)=\sum_{i=1}^{N_n}1_{[\#\{X_j\in\mathcal{P}_n|X_i\in
S_n(X_j,Y_j,t)\}\ge k_n+1]\bigcap[X_i\in A]} \nonumber
\end{equation}
Notice for the case $k_n\rightarrow\infty$, $s$ is suppressed in the
above expressions. Also, let
$\xi_n^{out}(t)\colon\hspace{-5pt}=\xi_n^{out}(t,\mathbb{R}^2)$
etc.\ for convenience.

The following two lemmas are intermediate steps to prove Theorem 1
and 2. We choose to state them without proof due to the limitation
of space and they can be treated in parallel with Theorem 4.12 and
4.13 in \cite{1} through a dependency graph argument.

\noindent\textbf{Lemma 1.} \itshape Suppose that $k_n=k$ is fixed,
and that $A$ is a Borel set in $\mathbb{R}^2$. The
finite-dimensional distributions of the process
$$n^{-\frac12}[\xi_n^{'out}(t,A)-E\xi_n^{'out}(t,A)]\quad , \quad t\ge 0$$
converge to those of a centered Gaussian process
$(\xi_{\infty}^{'out}(t,A),t>0)$ with covariance\\
$E[\xi_{\infty}^{'out}(t,A)\xi_{\infty}^{'out}(u,A)]$ given by
\begin{eqnarray}
\lefteqn{\int_A\rho_{\frac{\alpha}{2}
tf(x)}([k,\infty))f(x)\mathrm{d}x}\nonumber\\
&
&+\frac{1}{4\piup^2}\int_0^{2\piup}\int_0^{2\piup}\int_A\int_{\mathbb{R}^2}
\psi_{\infty}^{out}(z,f(x_1),y_1,y_2)f^2(x_1)\mathrm{d}z\mathrm{d}x_1\mathrm{d}y_1\mathrm{d}y_2\nonumber
\end{eqnarray}
with
\begin{eqnarray}
\psi_{\infty}^{out}(z,\lambda,y_1,y_2)&\hspace{-5pt}=&\hspace{-10pt}P(\{\mathcal{H}_{\lambda}^z(S(0,y_1,t^{\frac
12}))\ge k\}\cap\{\mathcal{H}_{\lambda}^0(S(z,y_2,u^{\frac
12}))\ge k\})\nonumber\\
&\hspace{-5pt}&\hspace{-10pt}-P(\mathcal{H}_{\lambda}(S(0,y_1,t^{\frac
12}))\ge k) P(\mathcal{H}_{\lambda}(S(z,y_2,u^{\frac 12}))\ge
k)\nonumber
\end{eqnarray}
The finite-dimensional distributions of
the process\\
$$n^{-\frac12}[\xi_n^{'in}(t,A)-E\xi_n^{'in}(t,A)]\quad , \quad t\ge 0.$$
converge to those of a centered Gaussian process
$(\xi_{\infty}^{'in}(t,A),t>0)$ with covariance\\
$E[\xi_{\infty}^{'in}(t,A)\xi_{\infty}^{'in}(u,A)]$ given by
\begin{equation}
\int_A\rho_{\frac{\alpha}{2}
tf(x)}([k,\infty))f(x)\mathrm{d}x+\int_A\int_{\mathbb{R}^2}
\psi_{\infty}^{in}(z,\frac{\alpha}{2\piup}f(x_1))f^2(x_1)\mathrm{d}z\mathrm{d}x_1\nonumber
\end{equation}
with
\begin{eqnarray}
\psi_{\infty}^{in}(z,\lambda)&\hspace{-5pt}=&\hspace{-10pt}P(\{\mathcal{H}_{\lambda}^z(B(0,t^{\frac
12}))\ge k\}\cap\{\mathcal{H}_{\lambda}^0(B(z,u^{\frac
12}))\ge k\})\nonumber\\
&\hspace{-5pt}&\hspace{-10pt}-P(\mathcal{H}_{\lambda}(B(0,t^{\frac
12}))\ge k) P(\mathcal{H}_{\lambda}(B(z,u^{\frac 12}))\ge
k)\nonumber
\end{eqnarray}\normalfont

Let $\mathcal{W}$ denote homogeneous white noise of intensity
$\piup^{-1}$ on $\mathbb{R}^2$, that is, a centered Gaussian process
indexed by bounded Borel sets in $\mathbb{R}^2$, with covariance
$\mathrm{Cov}(\mathcal{W}(A),\mathcal{W}(B))=\frac{1}{\piup}|A\cap
B|$, where $|\cdot|$ as mentioned before is Lebesgue measure. Also,
let $\mathcal{W}'$ denote homogeneous white noise of intensity
$\frac{2}{\alpha}$.

\medskip
\noindent\textbf{Lemma 2.} \itshape Suppose that
$k_n\rightarrow\infty$, that (\ref{1}) holds, and that $A$ is a
Borel set in $\mathbb{R}^2$. Let $s>0$ and suppose $F(A\cap L_s)>0$.
The finite-dimensional distributions of the process
$$(nk_n)^{-\frac12}[\xi_n^{'out}(t,A)-E\xi_n^{'out}(t,A)]\quad , \quad t\in\mathbb{R}$$
converge to those of a centered Gaussian process
$(\xi_{\infty}^{'out}(t,A),t\in\mathbb{R})$ with covariance
$E[\xi_{\infty}^{'out}(t,A)\xi_{\infty}^{'out}(u,A)]$ given by
\begin{equation}
\frac{|L_s\cap
A|}{s(\piup\alpha)^2}\int_0^{2\piup}\int_0^{2\piup}\int_{\mathbb{R}^2}\mathrm{Cov}
(1_{[\mathcal{W}'(S(0,y_1,1))\le t]},1_{[\mathcal{W}'(S(z,y_2,1))\le
u]})\mathrm{d}z\mathrm{d}y_1\mathrm{d}y_2\nonumber
\end{equation}
The finite-dimensional distributions of the process
$$(nk_n)^{-\frac12}[\xi_n^{'in}(t,A)-E\xi_n^{'in}(t,A)]\quad , \quad t\in\mathbb{R}$$
converge to those of a centered Gaussian process
$(\xi_{\infty}^{'in}(t,A),t\in\mathbb{R})$ with covariance\\
$E[\xi_{\infty}^{'in}(t,A)\xi_{\infty}^{'in}(u,A)]$ given by
\begin{equation}
\frac{4\cdot|L_s\cap A|}{s\alpha^2}\int_{\mathbb{R}^2}\mathrm{Cov}
(1_{[\mathcal{W}(B(0,1))\le t]},1_{[\mathcal{W}(B(z,1))\le
u]})\mathrm{d}z. \nonumber
\end{equation}\normalfont

Now we are ready to state our main results.

\medskip
\noindent\textbf{Theorem 1.} \itshape Suppose that $k_n=k$ is fixed.
The finite-dimensional distributions of the process
$$n^{-\frac12}[\xi_n^{out}(t)-E\xi_n^{out}(t)]\quad , \quad t\ge 0$$
converge to those of a centered Gaussian process
$(\xi_{\infty}^{out}(t),t>0)$ with
\begin{equation}
E[\xi_{\infty}^{out}(t)\xi_{\infty}^{out}(u)]=E[\xi_{\infty}^{'out}(t)\xi_{\infty}^{'out}(u)]
-h(t)h(u),\nonumber
\end{equation}
where
\begin{equation}
h(t)=\int_{\mathbb{R}^2}\left\{\rho_{\frac{\alpha}{2}
tf(x)}(k-1)\frac{\alpha}{2}tf(x)+\rho_{\frac{\alpha}{2}
tf(x)}([k,\infty))\right\}f(x)\mathrm{d}x\label{13}
\end{equation}
The above result also holds in the case where the superscripts 'out'
are replaced by 'in' everywhere. \normalfont

\medskip
\noindent\textbf{Theorem 2.} \itshape Suppose that $k_n
\rightarrow\infty$, and (\ref{1}) holds. Let $s>0$ and suppose
$F(L_s)>0$. The finite-dimensional distributions of the process
$$(nk_n)^{-\frac12}[\xi_n^{out}(t)-E\xi_n^{out}(t)]\quad , \quad t\in\mathbb{R}$$
converge to those of a centered Gaussian process
$(\xi_{\infty}^{out}(t),t\in\mathbb{R})$ with
\begin{equation}
E[\xi_{\infty}^{out}(t)\xi_{\infty}^{out}(u)]=E[\xi_{\infty}^{'out}(t)\xi_{\infty}^{'out}(u)]
-g(t)g(u),\nonumber
\end{equation}
where $g(t)=\phi(t)F(L_s)$.\\
The above result also holds in the case where the superscripts 'out'
are replaced by 'in' everywhere.\normalfont

\medskip
To deal with the degree distribution, let $\eta_n^{out}(t,A)$ and
$\eta_n^{in}(t,A)$ be the number of vertices in $A$ of out-degree
and in-degree $k$ fixed in $G_n$ respectively.

\medskip \noindent\textbf{Theorem 3.} \itshape Suppose $A$ is a Borel
set in $\mathbb{R}^2$ and $\alpha\ge\piup$. If either $k_n=k$ fixed,
or $k_n\rightarrow\infty$ and $n^{-1}k_n^2\ln n\rightarrow0$, then
\begin{equation}
\lim_{n\rightarrow\infty}n^{-1}\xi_n^{out}(t,A)-E[n^{-1}\xi_n^{out}(t,A)]=0\qquad
a.e.\label{60}
\end{equation}
Moreover,
\begin{equation}
\lim_{n\rightarrow\infty}n^{-1}\eta_n^{out}(t,A)=\int_A\rho_{\frac{\alpha}{2}
tf(x)}(k)f(x)\mathrm{d}x\qquad a.e.\label{61}
\end{equation}
The above result also holds in the case where the superscripts 'out'
are replaced by 'in' everywhere.\normalfont \bigskip

\noindent{\large\textbf{2.1 Discussion of Theorem 3.}}
\medskip

We take expectation on both sides of (\ref{61}), and let
$p(k):=E\lim_{n\rightarrow\infty}n^{-1}\eta_n^{out}(t,\mathbb{R}^2)$,
so the out-/in-degree distribution of
$G_{\alpha}(\mathcal{X}_n,\mathcal{Y}_n,r_n(t))$, where
$nr_n(t)^2=t$, is
\begin{equation}
p(k)=\frac{(\frac{\alpha}{2}t)^k}{k!}\int_{\mathbb{R}^2}e^{-\frac{\alpha}{2}tf(x)}f(x)^{k+1}\mathrm{d}x,\qquad
k\in\mathbb{N}\cup\{0\}\label{62}
\end{equation}

If we take the uniform density function $f(x)=1_{[0,1]^2}(x)$ in
(\ref{62}), then we see that
$p(k)=e^{-\frac{\alpha}{2}t}(\frac{\alpha}{2}t)^k/k!$, $k\ge0$; that
is, the degree distribution is $Poi(\frac{\alpha}{2}t)$.

If we take the standard multivariate normal density function
$f(x):=f(x_1,x_2)=(1/{2\piup})e^{-(x_1^2+x_2^2)/2}$, then through
the polar coordinate transformation and integration by parts, we
obtain $p(k)=(4\piup/\alpha t)-e^{-\alpha
t/4\piup}\sum_{i=0}^{k}(\alpha t/4\piup)^{i-1}/i!$, $k\ge0$. It is
easy to see that $p(k)\rightarrow0$ as $k\rightarrow\infty$; and
furthermore, since $p(0)=(4\piup/\alpha t)(1-e^{-\alpha t/4\piup})$,
$p(0)\rightarrow1$ as $t\rightarrow0$ and $p(0)\rightarrow0$ as
$t\rightarrow\infty$. These observations allow us presumably adjust
the parameter $t$ to get different skew degree distributions
especially for small $k$. However, the degree distribution in
(\ref{62}) has a light tail in contrast to the power law
distributions \cite{25} because of the fast decay as $k$ tends to
infinity. To be precise, by (\ref{62}) and Stirling formula,
$$
p(k)\le\frac{(\frac{\alpha}{2}tf_{\max})^k}{k!}\int_{\mathbb{R}^2}f(x)\mathrm{d}x
=(1+o(1))\cdot\frac{(\alpha tef_{\max})^k}{(2k)^k\sqrt{2\piup k}}
 \ll k^{-\beta}
$$
for any $\beta>0$ as $k\rightarrow\infty$.

On the other hand, if  we want to find a suitable density function
$f$ for a given probability distribution $p(k)$ satisfying $p(k)\ge
0$ and $\sum_{k=0}^{\infty}p(k)=1$, then we simply solve the
equation (\ref{62}), which is the first kind nonlinear singular
Fredholm integral equation \cite{26}. However, only approximation
solutions of this kind of equations may be obtained by using
iterative methods and the existence of solution is not known in
general.

\bigskip
\bigskip

\noindent{\Large\textbf{3. Proof of means and degree distribution}}

\medskip
\noindent\textbf{Proposition 1.} \itshape (out-degree) Suppose
$A\subseteq \mathbb{R}^2$ is a Borel set. If $k_n=k$ is fixed, then
\begin{equation}
\lim_{n\rightarrow
\infty}n^{-1}E[\xi_n^{out}(t,A)]=\int_A\rho_{\frac{\alpha}{2}
tf(x)}([k,\infty))f(x)\mathrm{d}x\label{15}
\end{equation}
If $k_n\rightarrow\infty$, and (\ref{1}) holds, then
\begin{equation}
\lim_{n\rightarrow \infty}n^{-1}E[\xi_n^{out}(t,A)]=F(L_s^+\cap
A)+\Phi(t)F(L_s\cap A)\label{16}
\end{equation}
\normalfont

\medskip
\noindent\textbf{Proof}. Let $p_n(x,y,t)=F(S_n(x,y,t))$. Then
\begin{equation}
E[\xi_n^{out}(t,A)]=\frac{n}{2\piup}\int_0^{2\piup}\int_A
P[Bin(n-1,p_n(x,y,t))\ge k_n]f(x)\mathrm{d}x\mathrm{d}y\label{17}
\end{equation}

Suppose $k_n$ is fixed, and $x\in R$, then $f$ is continuous at $x$
and $np_n(x,y,t)\rightarrow \frac{\alpha}{2}tf(x)$ by mean-value
theorem of integrals. Therefore $P[Bin(n-1,p_n(x,y,t))\ge k]$ tends
to $\rho_{\frac{\alpha}{2}tf(x)}([k,\infty))$. Then (\ref{15}) holds
by (\ref{17}) and dominated convergence theorem.

Suppose $k_n\rightarrow\infty$, (\ref{1}) holds and $x\in R$, then
$np_n(x,y,t)\sim n\frac{\alpha}{2}r_n^2f(x)\sim
s\frac{\alpha}{2}f(x)k_n$, and by Chernoff bounds (see
e.g.\cite{12}), $P[Bin(n-1,p_n(x,y,t))\ge k_n]$ tends to $1$, if $s
f(x)>\frac{2}{\alpha}$; and tends to $0$, if $s
f(x)<\frac{2}{\alpha}$. Then for $x\in R\cap L_s$,
\begin{eqnarray*}
np_n(x,y,t)&=&n\frac{\alpha}{2}
r_n^2f(x)+n\int_{S_n(x,y,t)}(f(z)-f(x))\mathrm{d}z\\
&=&k_n+tk_n^\frac 12+\Theta(n(k_n/n)^{3/2})
\end{eqnarray*}
Hence, by (\ref{1}),
\begin{equation}
np_n(x,y,t)=k_n+tk_n^\frac 12+o(k_n^{\frac 12}),\quad x\in R\cap
L_s\label{19}
\end{equation}

Then let $p_n=p_n(x,y,t)$, by DeMoivre-Laplace limit theorem and
(\ref{19}), we have
\begin{eqnarray*}
\lefteqn{P[Bin(n-1,p_n)\ge k_n]}\\
& &=P\left[\frac{Bin(n-1,p_n)-EBin(n-1,p_n)}{\sqrt{np_n}}\ge
\frac{k_n-(n-1)p_n}{\sqrt{np_n}}\right]\rightarrow \Phi(t)
\end{eqnarray*}
So (\ref{16}) follows from (\ref{17}) by dominated convergence
theorem. $\Box$

\medskip
\noindent\textbf{Proposition 2.} \itshape(in-degree) The same
results hold when replace superscripts ``out'' by ``in'' in
Proposition 1.\normalfont

\medskip
\noindent\textbf{Proof}. Let $q_n(x,t)=\frac{\alpha}{2\piup}\cdot
F(B_n(x,t))$. Then
\begin{equation}
E[\xi_n^{in}(t,A)]=n\int_A P[Bin(n-1,q_n(x,t))\ge
k_n]f(x)\mathrm{d}x.\nonumber
\end{equation}
From Palm theory, similarly we have
\begin{equation}
E[\xi_n^{'in}(t,A)]=n\int_A P[Poi(nq_n(x,t))\ge
k_n]f(x)\mathrm{d}x.\nonumber
\end{equation}
The remain proof is in a similar spirit with that of Proposition 1.
Hence we omit it. $\Box$
\medskip

We remark here that Proposition 1 and 2 still hold for corresponding
Poisson case.

\medskip
\noindent\textbf{Proof of Theorem 3}. Define a $\sigma$ filtration:
$\mathcal{F}_0=\{\emptyset,\Omega\}$, and for $1\le i\le n$,
$\mathcal{F}_i=\sigma\{(X_1,Y_1),(X_2,Y_2),\cdots,(X_i,Y_i)\}$.

For out-degree,
$\xi_n^{out}(t,A)-E[\xi_n^{out}(t,A)]=\sum_{i=1}^nM_{i,n}^{out}$,
with
$M_{i,n}^{out}=E[\xi_n^{out}(t,A)|\mathcal{F}_i]-E[\xi_n^{out}(t,A)|\mathcal{F}_{i-1}]$.
Let $\xi_{n,i}^{out}(t,A)$ be the number of vertices in $A$ of
$G(\mathcal{X}_{n+1}\backslash\{X_i\},\mathcal{Y}_{n+1}\backslash\{Y_i\}\\,r_n)$
having out-degree at least $k_n$. Thereby,
$M_{i,n}^{out}=E[\xi_n^{out}(t,A)-\xi_{n,i}^{out}(t,A)|\mathcal{F}_i]$.

We now claim that: For finite set $\mathcal{X}\subseteq\mathbb{R}^2$
and $x\in\mathcal{X}$, there are at most $8k$ points
$z\in\mathcal{X}$ having $x$ as their $(\le k)-th$ nearest neighbor,
for any $k\in\mathbb{N}$. Here $x$ is the $k-th$ nearest neighbor of
$z$ in $\mathcal{X}$ means if we order quantities $\{||w-z||:
w\in\mathcal{X}\backslash\{z\}\}$ increasingly, then $||x-z||$ will
be the $k-th$ item in this sequence. Proof. We take a cone with
vertex $x$, central angle $\piup/4$. It's easy to see that there are
at most $k_n$ points of $\mathcal{X}$ having $x$ as their $(\le
k)-th$ nearest neighbor, since we may look for these points from
near to far. The claim follows since the plane is covered by 8 such
cones.

Therefore,
\begin{eqnarray*}
|\xi_n^{out}(t,A)-\xi_{n,i}^{out}(t,A)|&\le&|\xi_n^{out}(t,A)-\tilde{\xi}_{n+1}^{out}(t,A)|+|\tilde{\xi}_{n+1}^{out}(t,A)-\xi_{n,i}^{out}(t,A)|\\
 &\le&(8k_n+1)+(8k_n+1)\le18k_n,
\end{eqnarray*}
where let $\tilde{\xi}_{n+1}^{out}(t,A)$ denote the number of
vertices in $A$ of out-degrees at least $k_n$ of
$G(\mathcal{X}_{n+1},\mathcal{Y}_{n+1},r_n)$. Then
$|M_{i,n}^{out}|\le 18k_n$. For $\varepsilon>0$, by Azuma
inequality, see e.g.\cite{24},
$$
P[|\xi_n^{out}(t,A)-E[\xi_n^{out}(t,A)]|>\varepsilon
n]\le2e^{-\varepsilon^2n^2/648nk_n^2}.
$$
By Borel-Cantelli Lemma, (\ref{60}) follows. The in-degree case can
be proved similarly.

To prove (\ref{61}), we notice
$$
\eta_n^{out}(t,A)=\sum_{i=1}^n1_{[\mathcal{X}_n(S_n(X_i,Y_i,t))\ge
k_n+1]\bigcap[X_i\in
A]}-\sum_{i=1}^n1_{[\mathcal{X}_n(S_n(X_i,Y_i,t))\ge
k_n+2]\bigcap[X_i\in A]}
$$
and by (\ref{60}) and the proof of Proposition 1, the result follows
immediately. The in-degree case also follows similarly. $\Box$

\bigskip
\bigskip
\noindent{\Large\textbf{4. Some moments for de-Poissonization}}
\smallskip

In this section we will develop some moments for non-Poisson case in
the limit regime $k_n\rightarrow\infty$, which is crucial to
de-Poisson Lemma 1 and 2.

For $n,m\in\mathbb{N}$, set
$$T_{m,n}^{out}(t):=\sum_{i=1}^{m}1_{[\mathcal{X}_m(S_n(X_i,Y_i,t)\backslash{\{X_i\}})\ge k_n]}$$
and
$$T_{m,n}^{in}(t):=\sum_{i=1}^{m}1_{[\#\{X_j\in\mathcal{X}_m\backslash{\{X_i\}}|X_i\in S_n(X_j,Y_j,t)\}\ge k_n]}$$
Then we see $T_{n,n}^{out}(t)=\xi_n^{out}(t),
T_{N_n,n}^{out}(t)=\xi_n^{'out}(t)$ and
$T_{n,n}^{in}(t)=\xi_n^{in}(t), T_{N_n,n}^{in}(t)=\xi_n^{'in}(t)$.
Set $\tilde{D}_{m,n}^{out}(t):=T_{m+1,n}^{out}(t)-T_{m,n}^{out}(t)$,
then
$\tilde{D}_{m,n}^{out}(t)=D_{m,n}^{out}(t)+\hat{D}_{m,n}^{out}(t)$,
where
$$D_{m,n}^{out}(t)=\sum_{i=1}^{m}1_{[\mathcal{X}_m(S_n(X_i,Y_i,t)\backslash{\{X_i\}})=k_n-1]\cap[X_{m+1}\in S_n(X_i,Y_i,t)]}$$
$$\hat{D}_{m,n}^{out}(t)=1_{[\mathcal{X}_m(S_n(X_{m+1},Y_{m+1},t))\ge k_n]}$$
Set $\tilde{D}_{m,n}^{in}(t):=T_{m+1,n}^{in}(t)-T_{m,n}^{in}(t)$,
then
$\tilde{D}_{m,n}^{in}(t)=D_{m,n}^{in}(t)+\hat{D}_{m,n}^{in}(t)$,
where
$$D_{m,n}^{in}(t)=\sum_{i=1}^{m}1_{[\#\{X_j\in\mathcal{X}_m\backslash{\{X_i\}}|X_i\in S_n(X_j,Y_j,t)\}=k_n-1]\cap[X_i\in S_n(X_{m+1},Y_{m+1},t)]}$$
$$\hat{D}_{m,n}^{in}(t)=1_{[\#\{X_j\in\mathcal{X}_m|X_{m+1}\in S_n(X_j,Y_j,t)\}\ge k_n]}$$

We denote binomial probability $\beta_{n,p}(k):=P(Bin(n,p)=k)$. The
next lemma will be repeatedly used in this section, see \cite{1,2}.

\medskip
\noindent\textbf{Lemma 3.} \itshape(a) Suppose $n,k\in \mathbb{N}$
with $k<n$. Then $\beta_{n,p}(k)$ is maximized over $p\in(0,1)$ by
setting $p=k/n$, and $p\beta_{n,p}(k)$ is maximized over $p\in(0,1)$
by setting $p=(k+1)/(n+1)$.

(b) Suppose $\{j_n\}_{n\ge1}$ is a sequence of integers satisfying
$j_n\rightarrow\infty$ and $(j_n/n)\rightarrow 0$ as
$n\rightarrow\infty$. Suppose $t\in\mathbb{R}$ and $\{p_n\}_{n\ge1}$
is a sequence in $(0,1)$ satisfying
$(j_n-np_n)/(np_n)^{1/2}\rightarrow t$ as $n\rightarrow\infty$. Then
$$j_n^{1/2}\beta_{n,p_n}(j_n)\rightarrow\phi(t) \qquad as\
n\rightarrow\infty.$$\normalfont

\medskip
\noindent\textbf{Lemma 4.} \itshape Suppose $k_n\rightarrow\infty$
and (\ref{1}) holds. Then
$$\lim_{n\rightarrow\infty}\sup_{\{m||m-n|\le n^{2/3}\}}|k_n^{-1/2}E\tilde{D}_{m,n}^{out}(t)-\phi(t)F(L_s)|=0$$
The same formula holds when replace superscript ``out'' by
``in''.\normalfont

\medskip
\noindent\textbf{Proof}. Take $\{m_n\}_{n\ge1}$ with $|m_n-n|\le
n^{2/3}$.

For out-degree, we have
\begin{eqnarray}
k_n^{-1/2}ED_{m_n,n}^{out}(t)&\hspace{-5pt}=&\hspace{-5pt}\frac{1}{2\piup}\int_0^{2\piup}\int_{\mathbb{R}^2}m_nk_n^{-1/2}P(\mathcal{X}_{m_n-1}(S_n(x,y,t))=k_n-1)\nonumber\\
& &\hspace{-5pt}\cdot
F(S_n(x,y,t))F(\mathrm{d}x)\mathrm{d}y\label{24}.
\end{eqnarray}
Let $x\in R\cap L_s$, then $\mathcal{X}_{m_n-1}(S_n(x,y,t))$ is
binomial with parameters $m_n-1$ and $F(S_n(x,y,t))$, and by
(\ref{19}), (\ref{1}) the mean is
\begin{eqnarray}
m_nF(S_n(x,y,t))&\hspace{-5pt}=&\hspace{-5pt}(1+O(n^{-1/3}))(k_n+tk_n^{1/2}+o(k_n^{1/2}))\nonumber\\
&\hspace{-5pt}=&\hspace{-5pt}k_n+tk_n^{1/2}+o(k_n^{1/2}),\qquad x\in
R\cap L_s\label{25}
\end{eqnarray}
By Lemma 3,
$$\lim_{n\rightarrow\infty}k_n^{1/2}P(\mathcal{X}_{m_n-1}(S_n(x,y,t))=k_n-1))=\phi(t),\quad
x\in R\cap L_s$$ Also, by Chernoff bounds and Proposition 1,
$$\lim_{n\rightarrow\infty}k_n^{1/2}P(\mathcal{X}_{m_n-1}(S_n(x,y,t))=k_n-1))=0,\quad
x\in R\backslash L_s$$ Hence for $x\in R$, the integrand on the
right hand side of (\ref{24}) tends to
$\frac{1}{2\piup}\phi(t)1_{L_s}(x)$. Also, by Lemma 3,
$(m_n/k_n)F(S_n(x,y,t))$ and
$k_n^{1/2}\sup_{0<p<1}\beta_{m_n-1,p}(k_n-1)$ are uniformly bounded.
So, $k_n^{-1/2}ED_{m_n,n}^{out}(t)$ tends to $\phi(t)F(L_s)$ by
dominated convergence theorem. Since $0\le
\hat{D}_{m_n,n}^{out}(t)\le 1$,
$k_n^{-1/2}E\hat{D}_{m_n,n}^{out}(t)$ tends to 0. The first part of
the lemma then follows.

For in-degree, we first introduce some notations. Let
$\tilde{f}:=\frac{\alpha}{2\piup}f$, and for Borel set $A\subseteq
\mathbb{R}^2$, let $\tilde{\mathcal{X}}_n(A)\sim
Bin(n,\tilde{F}(A))$, where
$\tilde{F}(A):=\int_A\tilde{f}(x)\mathrm{d}x$.

Consequently, we have
\begin{equation}
k_n^{-1/2}ED_{m_n,n}^{in}(t)=\int_{\mathbb{R}^2}m_nk_n^{-1/2}P(\tilde{\mathcal{X}}_{m_n-1}(B_n(x,t))=k_n-1)
\tilde{F}(B_n(x,t))F(\mathrm{d}x)\nonumber.
\end{equation}
Let $x\in R\cap L_s$, as mentioned above,
$\tilde{\mathcal{X}}_{m_n-1}(B_n(x,t))$ is binomial with parameters
$m_n-1$ and $\tilde{F}(B_n(x,t))$, and by Proposition 2 and
(\ref{1}) the mean is
\begin{equation}
m_n\tilde{F}(B_n(x,t))=k_n+tk_n^{1/2}+o(k_n^{1/2}),\qquad x\in R\cap
L_s\nonumber
\end{equation}

By using Lemma 3 and Proposition 2, we can conclude the proof in a
similar manner with the out-degree case. $\Box$

\medskip
\noindent\textbf{Lemma 5.} \itshape Suppose $k_n\rightarrow\infty$
and (\ref{1}) holds. Then
$$\lim_{n\rightarrow\infty}\sup_{n-n^{2/3}\le l<m\le n+n^{2/3}}|k_n^{-1}E\tilde{D}_{l,n}^{out}(t)\tilde{D}_{m,n}^{out}(u)-\phi(t)\phi(u)F(L_s)^2|=0$$
The same formula holds when replace superscripts ``out'' by
``in''.\normalfont

\medskip
\noindent\textbf{Proof}. Let $l\le m$.

For out degree, we have
\begin{eqnarray}
ED_{l,n}^{out}(t)D_{m,n}^{out}(u)&\hspace{-6pt}=&\hspace{-6pt}\sum_{i=1}^l\sum_{j=1}^mP[\{\mathcal{X}_l(S_n(X_i,Y_i,t))=k_n\}\cap\{\mathcal{X}_m(S_n(X_j,Y_j,u))=k_n\}\nonumber\\
& &\hspace{-6pt}\cap\{X_{l+1}\in S_n(X_i,Y_i,t)\}\cap\{X_{m+1}\in
S_n(X_j,Y_j,u)\}]\nonumber\\
&\hspace{-6pt}=&\hspace{-6pt}\frac{l(l-1)}{4\piup^2}\int_0^{2\piup}\int_0^{2\piup}\int_{\mathbb{R}^2}\int_{\mathbb{R}^2}g_{n,l,m}(x_1,y_1,x_2,y_2)F(\mathrm{d}x_1)F(\mathrm{d}x_2)\mathrm{d}y_1\mathrm{d}y_2\nonumber\\
& &\hspace{-6pt}+\frac{l(m-l)}{4\piup^2}\int_0^{2\piup}\int_0^{2\piup}\int_{\mathbb{R}^2}\int_{\mathbb{R}^2}g'_{n,l,m}(x_1,y_1,x_2,y_2)F(\mathrm{d}x_1)F(\mathrm{d}x_2)\mathrm{d}y_1\mathrm{d}y_2\nonumber\\
&
&\hspace{-6pt}+\frac{l}{2\piup}\int_0^{2\piup}\int_{\mathbb{R}^2}g''_{n,l,m}(x_1,y_1)F(\mathrm{d}x_1)\mathrm{d}y_1\label{28}
\end{eqnarray}
where,
\begin{eqnarray*}
g_{n,l,m}(x_1,y_1,x_2,y_2)&\hspace{-6pt}:=&\hspace{-6pt}P[\{\mathcal{X}_{l-2}^{x_2}(S_n(x_1,y_1,t))=k_n-1\}\cap\{X_{l-1}\in S_n(x_1,y_1,t)\}\\
&
&\hspace{-6pt}\cap\{\mathcal{X}_{m-2}^{x_1}(S_n(x_2,y_2,u))=k_n-1\}\cap\{X_{m-1}\in
S_n(x_2,y_2,u)\}]
\end{eqnarray*}
\begin{eqnarray*}
g'_{n,l,m}(x_1,y_1,x_2,y_2)&\hspace{-6pt}:=&\hspace{-6pt}P[\{\mathcal{X}_{l-1}(S_n(x_1,y_1,t))=k_n-1\}\cap\{X_{l}\in S_n(x_1,y_1,t)\}\\
&
&\hspace{-6pt}\cap\{\mathcal{X}_{m-2}^{x_1}(S_n(x_2,y_2,u))=k_n-1\}\cap\{X_{m-1}\in
S_n(x_2,y_2,u)\}]
\end{eqnarray*}
\begin{eqnarray*}
g''_{n,l,m}(x_1,y_1)&\hspace{-6pt}:=&\hspace{-6pt}P[\{\mathcal{X}_{l-1}(S_n(x_1,y_1,t))=k_n-1\}\cap\{X_{l}\in S_n(x_1,y_1,t)\}\\
&
&\hspace{-6pt}\cap\{\mathcal{X}_{m-1}(S_n(x_1,y_1,u))=k_n-1\}\cap\{X_{m}\in
S_n(x_1,y_1,u)\}]
\end{eqnarray*}

Take $x_1,x_2\in R$, $x_1\not=x_2$ and $y_1,y_2\in[0,2\piup)$. Take
$\{l_n\}_{n\ge1}$ and $\{m_n\}_{n\ge1}$ with $n-n^{2/3}\le
l_n<m_n\le n+n^{2/3}$. Then as $n\rightarrow\infty$,
\begin{equation}
\frac{n}{k_n}P(X_{l_n-1}\in
S_n(x_1,y_1,t))\rightarrow\frac{s\alpha}{2}f(x_1)\label{29},
\end{equation}
\begin{equation}
\frac{n}{k_n}P(X_{m_n-1}\in
S_n(x_2,y_2,u))\rightarrow\frac{s\alpha}{2}f(x_2)\label{30}.
\end{equation}
Since
\begin{multline}
P(\mathcal{X}_{m_n-2}^{x_1}(S_n(x_2,y_2,u))=k_n-1|X_{l_n-1}\in
S_n(x_1,y_1,t))\\
\sim\frac{P(\mathcal{X}_{m_n-2}(S_n(x_2,y_2,u))=k_n)\cdot
P(X_{l_n-1}\in S_n(x_1,y_1,t))}{P(X_{l_n-1}\in
S_n(x_1,y_1,t))}\sim\beta_{m_n,F(S_n(x_2,y_2,u))}(k_n)\nonumber,
\end{multline}
by Lemma 3 and (\ref{25}), we obtain
\begin{equation}
k_n^{1/2}P(\mathcal{X}_{m_n-2}^{x_1}(S_n(x_2,y_2,u))=k_n-1|X_{l_n-1}\in
S_n(x_1,y_1,t))\rightarrow\phi(u),\quad x_2\in L_s\label{31}
\end{equation}

Let $x_1, x_2\in\mathbb{R}^2$ with $x_2\not\in
B(x_1,r_n(t)+r_n(u))$, so $B_n(x_1,t)\cap B_n(x_2,u)=\emptyset$. If
$\mathcal{X}_{m_n-2}(S_n(x_2,y_2,u))=k_n-1$, then
$\mathcal{X}_{l_n-2}(S_n(x_2,y_2,u))=j$, for some $0\le j\le k_n-1$.
Given $\mathcal{X}_{l_n-2}(S_n(x_2,y_2,u))=j$, the conditional
distribution of $\mathcal{X}_{l_n-2}(S_n(x_1,y_1,t))$ is binomial
with parameter $l_n-2-j$ and
$F(S_n(x_1,y_1,t))/(1-F(S_n(x_2,y_2,u)))$. For all such $j$, if also
$x_1\in L_s$ then by (\ref{25}), (\ref{1}) the mean of this
distribution is
\begin{eqnarray*}
\frac{(l_n-2-j)F(S_n(x_1,y_1,t))}{1-F(S_n(x_2,y_2,u))}&\hspace{-6pt}=&\hspace{-6pt}(k_n+tk_n^{1/2}+o(k_n^{1/2}))(1+O(\frac{k_n}{n}))\\
&\hspace{-6pt}=&\hspace{-6pt}k_n+tk_n^{1/2}+o(k_n^{1/2})
\end{eqnarray*}
Therefore for $x_1\in L_s$ and $x_2\not=x_1$, by Lemma 3 we have
\begin{multline}
k_n^{1/2}P[\mathcal{X}_{l_n-2}^{x_2}(S_n(x_1,y_1,t))=k_n-1|\{\mathcal{X}_{m_n-2}^{x_1}(S_n(x_2,y_2,u))=k_n-1\}\\
\cap\{X_{l_n-1}\in S_n(x_1,y_1,t)\}]\rightarrow\phi(t).\nonumber
\end{multline}
Combining this with (\ref{29}), (\ref{30}) and (\ref{31}), we get
\begin{equation}
(n^2/k_n)g_{n,l_n,m_n}(x_1,y_1,x_2,y_2)\rightarrow\phi(t)\phi(u),\quad
x_1, x_2\in R\cap L_s, x_1\not=x_2\label{32}
\end{equation}
On the other hand, by Chernoff bounds,
$$k_nP[\mathcal{X}_{m_n-2}^{x_1}(S_n(x_2,y_2,u))=k_n-1|X_{l_n-1}\in S_n(x_1,y_1,t)]\rightarrow 0,\quad x_2\in R\backslash L_s$$
and
$$k_nP[\mathcal{X}_{l_n-2}^{x_2}(S_n(x_1,y_1,t))=k_n-1]\rightarrow 0,\quad x_1\in R\backslash L_s$$
Combing these with (\ref{29}) and (\ref{30}), we have
\begin{equation}
(n^2/k_n)g_{n,l_n,m_n}(x_1,y_1,x_2,y_2)\rightarrow 0, \quad (x_1,
x_2)\in (R\times R)\backslash (L_s\times L_s)\label{33}
\end{equation}

If $x_2\not\in B(x_1,r_n(t)+r_n(u))$, setting
$p_1=F(S_n(x_2,y_2,u))$ and $p_2=F(S_n(x_1,y_1,t))/(1-p_1)$, we get
$$g_{n,l_n,m_n}(x_1,y_1,x_2,y_2)\le \max_{0\le j\le k_n-1}p_1p_2\beta_{m_n-2,p_1}(k_n-1)\beta_{l_n-2-j,p_2}(k_n-1)$$
Whence by Lemma 3 and Stirling formula, there exists a constant $c$
such that
$$g_{n,l_n,m_n}(x_1,y_1,x_2,y_2)\le c(k_n^{-\frac12}\cdot\frac{k_n}{n})^2=ck_n/n^2.$$
Then by (\ref{32}), (\ref{33}) and dominated convergence theorem, we
get
\begin{multline}
\lim_{n\rightarrow\infty}\big[\frac{n^2}{k_n\cdot4\piup^2}\int_0^{2\piup}\int_0^{2\piup}\int_{\mathbb{R}^2}\int_{\mathbb{R}^2\backslash
B(x_1,r_n(t)+r_n(u))}g_{n,l_n,m_n}(x_1,y_1,x_2,y_2)\\
F(\mathrm{d}x_2)F(\mathrm{d}x_1)\mathrm{d}y_1\mathrm{d}y_2\big]=\phi(t)\phi(u)F(L_s)^2\label{34}.
\end{multline}
Also, by (\ref{29}), (\ref{30}) and (\ref{31}),
\begin{eqnarray*}
g_{n,l_n,m_n}(x_1,y_1,x_2,y_2)&\hspace{-6pt}\le&\hspace{-6pt}P(\mathcal{X}_{l_n-2}^{x_2}(S_n(x_1,y_1,t))=k_n-1)F(S_n(x_1,y_1,t))F(S_n(x_2,y_2,u))\\
&\hspace{-6pt}=&\hspace{-6pt}O\big(k_n^{-\frac12}\cdot\big(\frac{k_n}{n}\big)^2\big)
\end{eqnarray*}
Since $F(B(x_1,r_n(t)+r_n(u)))\le c\cdot (k_n/n)$ for some constant
$c$, by (\ref{1}),
\begin{multline}
\big(\frac{n^2}{k_n\cdot4\piup^2}\big)\int_0^{2\piup}\int_0^{2\piup}\int_{\mathbb{R}^2}\int_{B(x_1,r_n(t)+r_n(u))}g_{n,l_n,m_n}(x_1,y_1,x_2,y_2)\\
F(\mathrm{d}x_2)F(\mathrm{d}x_1)\mathrm{d}y_1\mathrm{d}y_2\le
c'\big(\frac{n^2}{k_n}\big)\big(\frac{k_n}{n}\big)\big(\frac{k_n^{3/2}}{n^2}\big)\rightarrow
0\nonumber
\end{multline}
Thus (\ref{34}) holds with the region of integration modified to
$[0,2\piup)\times[0,2\piup)\times\mathbb{R}^2\times\mathbb{R}^2$.
The asymptotic results for $g'_{n,l_n,m_n}$ are just the same. Also,
by similar arguments there is a constant $c$ such that
$$l_n\cdot\sup_{x_1\in R \atop y_1\in [0,2\piup)}g''_{n,l_n,m_n}(x_1,y_1)\le cnk_n^{-1/2}(k_n/n)^2\rightarrow 0.$$
Hence (\ref{28}) yields
$$k_n^{-1}ED_{l_n,n}^{out}(t)D_{m_n,n}^{out}(u)\rightarrow \phi(t)\phi(u)F(L_s)^2.$$
What remains to show is that the above formula still holds when
$D_{l_n,n}^{out}$ is replaced by $\tilde{D}_{l_n,n}^{out}$;
$D_{m_n,n}^{out}$ is replaced by $\tilde{D}_{m_n,n}^{out}$. We have
$0\le \hat{D}_{l_n,n}^{out}(t)\le 1$, $0\le
\hat{D}_{m_n,n}^{out}(u)\le 1$. By the proof of Lemma 4,
$ED_{l_n,n}^{out}(t)=O(k_n^{1/2})$ and
$ED_{m_n,n}^{out}(u)=O(k_n^{1/2})$. Therefore
$E[D_{l_n,n}^{out}(t)\hat{D}_{m_n,n}^{out}(u)]$,
$E[\hat{D}_{l_n,n}^{out}(t)\hat{D}_{m_n,n}^{out}(u)]$ and
$E[\hat{D}_{l_n,n}^{out}(t)D_{m_n,n}^{out}(u)]$ are all
$O(k_n^{1/2})$. The first part of this lemma whereby follows.

The proof for in-degree case parallels to the above approach and we
leave it as an exercise for the reader.$\Box$

\medskip
\noindent\textbf{Lemma 6.} \itshape Suppose $k_n\rightarrow\infty$
and (\ref{1}) holds. Let $t,u\in\mathbb{R}$. Then
$$\limsup_{n\longrightarrow\infty}\big(k_n^{-3/2}\cdot\sup_{\{m||m-n|\le n^{2/3}\}}E[\tilde{D}^{out}_{m,n}(t)^2]\big)<\infty.$$
The same formula holds when replace superscript ``out'' by
``in''.\normalfont

\medskip
\noindent\textbf{Proof}. Take $\{m_n\}_{n\ge1}$ satisfying
$|m_n-n|\le n^{2/3}$.

For out-degree, by (\ref{28}) with $l=m=m_n$, $t=u$,
\begin{eqnarray*}
E[D_{m_n,n}^{out}(t)^2]&\hspace{-6pt}=&
\hspace{-6pt}\frac{m_n(m_n-1)}{4\piup^2}\int_0^{2\piup}\int_0^{2\piup}\int_{\mathbb{R}^2}\int_{\mathbb{R}^2}g_{n,m_n,m_n}(x_1,y_1,x_2,y_2)F(\mathrm{d}x_1)F(\mathrm{d}x_2)\mathrm{d}y_1\mathrm{d}y_2\\
& &\hspace{-6pt}+ED_{m_n,n}^{out}(t)
\end{eqnarray*}
By Lemma 3, there is a constant $c$ such that
$$g_{n,m_n,m_n}(x_1,y_1,x_2,y_2)\le \sup_{0<p<1}[p\max\{\beta_{m_n,p}(k_n-1),\beta_{m_n,p}(k_n-2)\}]\le\big(\frac{ck_n}{n}\big)k_n^{-1/2}.$$
Also, $g_{n,m_n,m_n}(x_1,y_1,x_2,y_2)=0$ unless $x_2\in
B(x_1,2r_n(t))$. Whence
$$\int_0^{2\piup}\int_0^{2\piup}\int_{\mathbb{R}^2}\int_{\mathbb{R}^2}g_{n,m_n,m_n}(x_1,y_1,x_2,y_2)F(\mathrm{d}x_1)F(\mathrm{d}x_2)\mathrm{d}y_1\mathrm{d}y_2\le \frac{c'k_n^{3/2}}{n^2}$$
By Lemma 4, $ED_{m_n,n}^{out}(t)=O(k_n^{1/2})$ and
$m_n(m_n-1)=O(n^2)$. So $E[D_{m_n,n}^{out}(t)^2]=O(k_n^{3/2})$. The
first part of this lemma then follows, by noting $0\le
\hat{D}_{m_n,n}^{out}(t)\le 1$.

For in-degree, the same argument may be applied. Thus we conclude
the proof. $\Box$
\bigskip
\bigskip

\noindent{\Large\textbf{5. Proof of central limit theorems}}
\smallskip

\medskip
To prove Theorem 1 and 2, we will employ useful de-Poisson
techniques given in \cite{16}, \cite{17} and later generalized in
\cite{1,19}. We will also need Cram\'er-Wold device, see
e.g.\cite{12}. Now we are in position to prove our main results.

\medskip
\noindent\textbf{Proof of Theorem 1.} Let $M\in\mathbb{N}$,
$B=(b_1,\cdots,b_M)\in \mathbb{R}^M$, $T=(t_1,\cdots,t_M)\in
(0,\infty)^M$.

For out-degree,
$\mathcal{X}\subset\mathbb{R}^2,\mathcal{Y}\subset[0,2\piup)$ with
$\mathrm{card}(\mathcal{X})=\mathrm{card}(\mathcal{Y})$, set
$$H_0(\mathcal{X},\mathcal{Y}):=\sum_{i=1}^M\sum_{(x,y)\in(\mathcal{X},\mathcal{Y})}b_i1_{[\mathcal{X}(S(x,y,t_i^{1/2}))\ge
k_n+1]}$$ and let
$H_n(\mathcal{X},\mathcal{Y})=H_0(n^{1/2}\mathcal{X},\mathcal{Y})$.
$(x,y)\in\mathbb{R}^2\times[0,2\piup)\subset\mathbb{R}^{3}$. Set
$\xi_n^{'out}(T,B,A):=\sum_{m=1}^Mb_m\xi_n^{'out}(t_m,A)$ and
$\mathrm{Var}\big(\xi_n^{'out}(T,B,A)\big):=\sigma^{'out}(T,B,A)$,
we have
$H_n(\mathcal{P}_n,\mathcal{Y}_{N_n})=\xi_n^{'out}(T,B,\mathbb{R}^2)$,
and what's more, $(\mathcal{P}_n,\mathcal{Y}_{N_n})$ is a
$3-$dimensional Poisson process, which may be coupled with
$(\mathcal{X}_n,\mathcal{Y}_n)$ in the same way as $\mathcal{P}_n$
does with $\mathcal{X}_n$. By Lemma 1,
$n^{-1/2}(H_n(\mathcal{P}_n,\mathcal{Y}_{N_n})-EH_n(\mathcal{P}_n,\mathcal{Y}_{N_n}))\stackrel{\mathrm{D}}{\longrightarrow}\mathcal{N}(0,\sigma^{'out}(T,B,\mathbb{R}^2)).$
Let $\mathcal{H}_{\lambda}$ be a $3-$dimensional homogeneous Poisson
process and denote point $(x,y)\in\mathbb{R}^2\times\mathbb{R}$,
$\mathcal{H}_{\lambda}:=(\mathcal{H}^{(1)}_{\lambda},\mathcal{H}^{(2)}_{\lambda})$
with
$x\in\mathcal{H}^{(1)}_{\lambda},y\in\mathcal{H}^{(2)}_{\lambda}$.
Next, we say $H_0(\mathcal{X},\mathcal{Y})$ is strongly stabilizing
on $\mathcal{H}_{\lambda}$ if there are $a.s.$ finite random
variables $T$ and $\Delta(\mathcal{H}_{\lambda})$ such that with
probability 1, $\Delta(A)=\Delta(\mathcal{H}_{\lambda})$ for all
finite $A:=(A_1,A_2)\subset\mathbb{R}^2\times[0,2\piup)$ with
$\mathrm{card}(A_1)=\mathrm{card}(A_2)$, satisfying
$A\cap(B(0,T)\times[0,2\piup))=\mathcal{H}_{\lambda}\cap(B(0,T)\times[0,2\piup))$.
Here,
$\triangle(\mathcal{H}_{\lambda}):=H_0(\mathcal{H}^0_{\lambda})-H_0(\mathcal{H}_{\lambda})$.
Thus, $H_0$ is strongly stable since it has finite range. We have
\begin{eqnarray*}
E[\triangle(\mathcal{H}_{\lambda})]&\hspace{-6pt}=&\hspace{-6pt}E[H_0(\mathcal{H}^0_{\lambda})-H_0(\mathcal{H}_{\lambda})]\\
&\hspace{-6pt}=&\hspace{-6pt}E\big[\sum_{i=1}^Mb_i\big(\sum_{(x,y)\in\mathcal{H}_{\lambda}}1_{[\mathcal{H}^{(1),0}_{\lambda}(S(x,y,t_i^{1/2}))\ge
k+1]}+1_{[\mathcal{H}^{(1),0}_{\lambda}(S(0,0,t_i^{1/2}))\ge
k+1]}\big)\\
&\hspace{-6pt}
&\hspace{-6pt}-\sum_{i=1}^Mb_i\big(\sum_{(x,y)\in\mathcal{H}_{\lambda}}1_{[\mathcal{H}^{(1)}_{\lambda}(S(x,y,t_i^{1/2}))\ge
k+1]}\big)\big]\\
&\hspace{-6pt}=&\hspace{-6pt}E\sum_{i=1}^Mb_i\big(1_{[Poi(2\piup\lambda\cdot\frac{\alpha
t_i}{2})\ge k]}+\sum_{(x,y)\in\mathcal{H}_{\lambda} \atop 0\in
S(x,y,t_i^{1/2})}1_{[\mathcal{H}^{(1)}_{\lambda}(S(x,y,t_i^{1/2}))=
k]}\big)\\
&\hspace{-6pt}=&\hspace{-6pt}\sum_{i=1}^Mb_i\big(\rho_{\lambda\piup\alpha
t_i}([k,\infty))+\lambda2\piup\cdot\frac{\alpha t_i}{2}(k-1)\big).
\end{eqnarray*}
By (\ref{13}) and the Cox process $\mathcal{H}_{\varphi(X,Y)}$ with
$\varphi(X,Y):=(1/2\piup)f(X)$, we have
$E[\triangle(\mathcal{H}_{\varphi(X,Y)})]\\=\sum_{i=1}^Mb_ih(t_i)$.
Set $t_{\mathrm{max}}=\max\{t_1,\cdots,t_M\}$, we have
$|H_n(\mathcal{X}_m,\mathcal{Y}_m)|\le m\sum_{i=1}^M|b_i|$ and
$$|H_n(\mathcal{X}_{m+1},\mathcal{Y}_{m+1})-H_n(\mathcal{X}_m,\mathcal{Y}_m)|\le \big(\sum_{i=1}^Mb_i\big)\cdot\big[\#\{X_i\in\mathcal{X}_m|X_{m+1}\in S_n(X_i,Y_i,t_{\mathrm{max}})\}+1\big]$$
which is stochastically dominated by
$c\cdot[Bin(m,f_{\mathrm{max}}\piup r_n(t_{\mathrm{max}})^2)+1]$
having a uniformly bounded fourth moment when $m\le2n$. Therefore by
a simple variant of Theorem 2.16(\cite{1}) to a marked point process
\cite{19} (in particular the translation-invariance of
$H_0(\mathcal{X},\mathcal{Y})$ is only required for $\mathcal{X}$),
$n^{-1/2}(H_n(\mathcal{X}_n,\mathcal{Y}_n)-EH_n(\mathcal{X}_n,\mathcal{Y}_n))\stackrel{\mathrm{D}}{\longrightarrow}\mathcal{N}(0,\tau_{out}^2)$
with
$\tau_{out}^2:=\sigma^{'out}(T,B,\mathbb{R}^2)-(E[\triangle(\mathcal{H}_{\varphi(X,Y)})])^2$.
The first part of the theorem then follows by Cram\'er-Wold device.

For in-degree, let
$\mathcal{X}\subset\mathbb{R}^2,\mathcal{Y}\subset[0,2\piup)$ with
$\mathrm{card}(\mathcal{X})=\mathrm{card}(\mathcal{Y})$, and the
elements $(x,y)\in(\mathcal{X},\mathcal{Y})$ be ordered pairs. Reset
$$H_0(\mathcal{X},\mathcal{Y}):=\sum_{i=1}^M\sum_{x'\in\mathcal{X}}b_i1_{[\#\{x\in\mathcal{X}|x'\in S(x,y,t_i^{1/2})\}\ge
k_n+1]}$$ and let
$H_n(\mathcal{X},\mathcal{Y})=H_0(n^{1/2}\mathcal{X},\mathcal{Y})$.
Set $\xi_n^{'in}(T,B,A):=\sum_{m=1}^Mb_m\xi_n^{'in}(t_m,A)$ and\\
$\mathrm{Var}\big(\xi_n^{'in}(T,B,A)\big):=\sigma^{'in}(T,B,A)$, we
have
$H_n(\mathcal{P}_n,\mathcal{Y}_{N_n})=\xi_n^{'in}(T,B,\mathbb{R}^2)$,
and $(\mathcal{P}_n,\mathcal{Y}_{N_n})$ is a $3-$dimensional Poisson
process coupled with $(\mathcal{X}_n,\mathcal{Y}_n)$. By Lemma 1,
$n^{-1/2}(H_n(\mathcal{P}_n,\mathcal{Y}_{N_n})\\-EH_n(\mathcal{P}_n,\mathcal{Y}_{N_n}))\stackrel{\mathrm{D}}{\longrightarrow}\mathcal{N}(0,\sigma^{'in}(T,B,\mathbb{R}^2)).$
Also, $H_0$ is strongly stable. Let $\mathcal{H}_{\lambda}$ be a
$3-$dimensional homogeneous Poisson process and
$\mathcal{H}_{\lambda}:=(\mathcal{H}^{(1)}_{\lambda},\mathcal{H}^{(2)}_{\lambda})$
as above. Then
\begin{eqnarray*}
E[\triangle(\mathcal{H}_{\lambda})]&\hspace{-7pt}=&\hspace{-7pt}E[H_0(\mathcal{H}^0_{\lambda})-H_0(\mathcal{H}_{\lambda})]\\
&\hspace{-7pt}=&\hspace{-7pt}E\big[\sum_{i=1}^Mb_i\big(\sum_{x'\in\mathcal{H}^{(1)}_{\lambda}}1_{[\#\{x\in\mathcal{H}_{\lambda}^{(1),0}|x'\in
S(x,y,t_i^{1/2})\}\ge
k+1]}+1_{[\#\{x\in\mathcal{H}_{\lambda}^{(1),0}|0\in
S(x,y,t_i^{1/2})\}\ge k+1]}\big)\\
&\hspace{-7pt}
&\hspace{-7pt}-\sum_{i=1}^Mb_i\big(\sum_{x'\in\mathcal{H}^{(1)}_{\lambda}}1_{[\#\{x\in\mathcal{H}_{\lambda}^{(1)}|x'\in
S(x,y,t_i^{1/2})\}\ge
k+1]}\big)\big]\\
&\hspace{-7pt}=&\hspace{-7pt}E\sum_{i=1}^Mb_i\big(1_{[Poi(2\piup\lambda\cdot\frac{\alpha
t_i}{2})\ge k]}+\sum_{x'\in\mathcal{H}^{(1)}_{\lambda}\cap
S(0,0,t_i^{1/2})}1_{[\#\{x\in\mathcal{H}_{\lambda}^{(1)}|x'\in
S(x,y,t_i^{1/2})\}=k]}\big)\\
&\hspace{-7pt}=&\hspace{-7pt}\sum_{i=1}^Mb_i\big(\rho_{\lambda\piup\alpha
t_i}([k,\infty))+2\piup\lambda\cdot\frac{\alpha t_i}{2}(k-1)\big).
\end{eqnarray*}
The remain proof is similar with the out-degree case. $\Box$

\medskip
\noindent\textbf{Proof of Theorem 2.} Let $T$ and
$B\in\mathbb{R}^M$.

For out-degree,
$\mathcal{X}\subset\mathbb{R}^2,\mathcal{Y}\subset[0,2\piup)$ with
$\mathrm{card}(\mathcal{X})=\mathrm{card}(\mathcal{Y})$, set
$$H_n(\mathcal{X},\mathcal{Y}):=k_n^{-1/2}\sum_{i=1}^M\sum_{(x,y)\in(\mathcal{X},\mathcal{Y})}b_i1_{[\mathcal{X}(S_n(x,y,t_i))\ge k_n+1]}$$
By Lemma 2, we have
$H_n(\mathcal{P}_n,\mathcal{Y}_{N_n})=k_n^{-1/2}\xi_n^{'out}(T,B,\mathbb{R}^2)$
and
$n^{-1/2}(H_n(\mathcal{P}_n,\mathcal{Y}_{N_n})-EH_n(\mathcal{P}_n,\mathcal{Y}_{N_n}))\stackrel{\mathrm{D}}{\longrightarrow}\mathcal{N}(0,\sigma^{'out}(T,B,\mathbb{R}^2)).$
Set $\alpha:=\sum_{i=1}^Mb_i\phi(t_i)F(L_s)$, and
$R^{out}_{m,n}:=H_n(\mathcal{X}_{m+1},\mathcal{Y}_{m+1})-H_n(\mathcal{X}_m,\mathcal{Y}_m)$.
Then
$R^{out}_{m,n}=k_n^{-1/2}\sum_{i=1}^Mb_i\tilde{D}_{m,n}^{out}(t_i).$
By Lemma 4, 5 and 6, we have
$$\lim_{n\rightarrow\infty}\big(\sup_{n-n^{2/3}\le m\le n+n^{2/3}}|ER^{out}_{m,n}-\alpha|\big)=0$$
$$\lim_{n\rightarrow\infty}\big(\sup_{n-n^{2/3}\le m<m'\le n+n^{2/3}}|E[R^{out}_{m,n}R^{out}_{m',n}]-\alpha^2|\big)=0$$
and
$$\lim_{n\rightarrow\infty}\big(n^{-1/2}\sup_{n-n^{2/3}\le m\le n+n^{2/3}}E[(R^{out}_{m,n})^2]\big)=0.$$
respectively. Also $|H_n(\mathcal{X}_m,\mathcal{Y}_m)|\le
m\sum_{i=1}^M|b_i|$. Then Theorem2.12(\cite{1}) implies
$n^{-1/2}\\\cdot(H_n(\mathcal{X}_n,\mathcal{Y}_n)-EH_n(\mathcal{X}_n,\mathcal{Y}_n))\stackrel{\mathrm{D}}{\longrightarrow}\mathcal{N}(0,\sigma^{out}(T,B))$,
with $\sigma^{out}(T,B):=\sigma^{'out}(T,B,\mathbb{R}^2)-\alpha^2.$
Hence,
$\sigma^{out}(T,B)=\mathrm{Var}\sum_{i=1}^Mb_i\xi_{\infty}^{out}(t_i)$.
The first part of the theorem then follows by Cram\'er-Wold device.

For in-degree, let
$$H_n(\mathcal{X},\mathcal{Y}):=k_n^{-1/2}\sum_{i=1}^M\sum_{x'\in\mathcal{X}}b_i1_{[\#\{x\in\mathcal{X}|x'\in S_n(x,y,t_i)\}\ge
k_n+1]}$$ We then argue likewise to complete the proof. $\Box$

\bigskip
\bigskip

\noindent{\Large\textbf{6. Further discussion and remarks}}
\smallskip

In the above sections, we consider $d=2$ and $Y_i$ uniformly
distributed. A natural generalization is to consider higher
dimensions. For example, for $d=3$, instead of a sector with
amplitude $\alpha$, we have to consider a spherical sector
$SS(X,Y,Z,r)$ which is the region bounded by a cone with vertex $X$,
central angle $\alpha$ and a sphere with center $X$ and radius $r$.
We take $X$ as the origin and build the standard right-handed
coordinate system. Let the chief axis of the cone be a ray $l$,
project $l$ onto $xOy-$plane, and call it $l'$. Let $Y$ be the angle
between positive $x-$axis and $l'$ and $Z+(\alpha/2)$ be the angle
between $l$ and $l'$. $Y,Z\in[0,2\piup)$. Consequently, the formal
definition of this ``random spherical sector graph'' is easily
stated. If $Y$ and $Z$ have uniform distribution, and instead of
condition (\ref{1}) we assume $k_n/n^{\frac2{d+2}}$ tends to $0$ and
modify the definitions of $r_n(t)$ accordingly, then analogous
results corresponding to those appeared in above sections can be
derived. Actually, we have for example,
$\xi_n^{out}(t,A)=\sum_{i=1}^n1_{[\mathcal{X}_n(SS_n(X_i,Y_i,Z_i,t))\ge
k_n+1]\bigcap[X_i\in A]}$. Let $p_n(x,y,z,t)=F(SS_n(x,y,z,t))$, then
$$E[\xi_n^{out}(t,A)]=\frac{n}{4\piup^2}\int_0^{2\piup}\int_0^{2\piup}\int_A
P[Bin(n-1,p_n(x,y,z,t))\ge
k_n]f(x)\mathrm{d}x\mathrm{d}y\mathrm{d}z.$$ Also,
$\xi_n^{in}(t,A)=\sum_{i=1}^n1_{[\#\{X_j\in\mathcal{X}_n|X_i\in
SS_n(X_j,Y_j,Z_j,t)\}\ge k_n+1]\bigcap[X_i\in A]}$. Let
$q_n(x,t)=F(B_n(x,t))\\\cdot\big(\frac{1-\cos(\alpha/2)}{2}\big)$,
then
$$E[\xi_n^{in}(t,A)]=n\int_A P[Bin(n-1,q_n(x,t))\ge
k_n]f(x)\mathrm{d}x.$$

Another direction to investigate is to consider probability density
 $g$ of $Y$ other than the uniform density. Suppose $EY<\infty$. For
out-degree case, we may proceed smoothly by similar argument,
whereas for in-degree the story is different. Say, we consider
in-degree of a vertex $u$. Suppose $\|u-v\|<r.$ Since the
inclination of sector $S_v$ now is not uniformly at random (as we
now consider a general density $g$), we will have distinct thinning
probability for different $v$. Moreover, the probability of vertex
$u$ lying in the sector $S_v$ essentially relies on not only the
distance between them but also the position of both vertices $u$ and
$v$. Then the computation is inevitably involved and the above
de-Poisson technique is no longer valid.

We mention that the model is less interesting when using other
non-Euclidean norm in application viewpoint. It is easy to see when
$d=2$, if we take $l^p\ (1\le p\le\infty)$ norm, and
$\alpha=\piup/2,\piup,3\piup/2$ or $2\piup$, the above results still
hold, due to the symmetry of the coordinate vectors under such norm.

\bigskip

\noindent{\Large\textbf{Acknowledgements}}

\smallskip
I would like to thank Professor Qing Zhou and Professor Dong Han for
helpful comments and sound advice.


\bigskip

\begin{thebibliography}{10}
\bibitem{23}B.\ Bollob\'as, \textit{Random Graphs}. Cambridge
University Press, 2001
\bibitem{24}F.\ R.\ K.\ Chung, Linyuan Lu, \textit{Complex Graphs and
Networks}. American Mathematical Society, CBMS. 2006
\bibitem{26}P.\ J.\ Collins, \textit{Differential and Integral
Equations}. Oxford University Press, 2006
\bibitem{3}J.\ D\'{\i}az, J.\ Petit, M.\ Serna, A random graph model
for optical networks of sensors. \textit{IEEE Transactions on Mobile
Computing} 2(2003) 143--154
\bibitem{4}J.\ D\'{\i}az, V.\ Sanwalani, M.\ Serna, P.\ G.\
Spirakis, The chromatic and clique numbers of random scaled sector
graphs. \textit{Theoretical Computer Science} 349(2005) 40--51
\bibitem{18}J.\ D\'{\i}az, Zvi Lotker, M.\ Serna, The distant-2 chromatic number of random
proximity and random geometric graphs. \textit{Information
Processing Letters} 106(2008) 144--148
\bibitem{25}S.\ N.\ Dorogovstev, J.\ F.\ F.\ Mendes, Evolution of
networks, \textit{Advances in Physics} 51(2002) 1079--1187
\bibitem{12}S.\ Janson, T.\ {\L}uczak, A.\ Rucinski, \textit{Random
Graphs}. Wiley, New York, 2000
\bibitem{16}S.\ Lee, The central limit theorem for Euclidean minimal
spanning trees. I. \textit{Ann.\ Appl.\ Probab.\ }7(1997) 996--1020
\bibitem{5}T.\ M\H{u}ller, Two point concentration in random geometric
graphs, \textit{Combinatorica}, 2009 article in press, 1--17
\bibitem{20}S.\ Nikoletseas, Models and algorithms for wireless sensor networks (Smart
Dust), \textit{LNCS} 3831(2006) 64--83
\bibitem{7}T.\ Pei, A.-X. Zhu, C.\ Zhou, B.\ Li, C. Qin, A new
approach to the nearest-neighbour method to discover cluster
features in overlaid spatial point processes, \textit{International
Journal of Geographical Information Science} 20(2006) 153--168
\bibitem{1}M.\ D.\ Penrose, \textit{Random Geometric Graphs}. Oxford
University Press, Oxford, 2003
\bibitem{2}M.\ D.\ Penrose. Central limit theorems for $k$-nearest
neighbor distances. \textit{Stochastic Processes and their
Applications} 85(2000) 295--320
\bibitem{17}M.\ D.\ Penrose, J.\ E.\ Yukich, Central limit theorems
for some graphs in computational geometry. \textit{Ann.\ Appl.\
Probab.\ } 11(2001) 1005--1041
\bibitem{19}M.\ D.\ Penrose, Multivariate spatial central limit
theorems with applications to percolation and spatial graphs.
\textit{Annals of Probability} 33(2005) 1945--1991
\bibitem{22}O.\ Schabenberger, C.\ A.\ Gotway, \textit{Statistical Methods for Spatial Data
Analysis}. CRC Press, 2004
\bibitem{21}M.\ Serna, Random models for geometric graphs, \textit{LNCS} 4525,
pp.37, 2007
\bibitem{30}Y.\ Shang, Focusing of maximum vertex degrees in random faulty scaled sector
graphs, arXiv:0909.2933v1 [math.CO].
\end{thebibliography}
\end{document}